\documentstyle[12pt,amsfonts]{article}

\begin{document}
\title{Regularity of  H\"{o}lder continuous solutions of the supercritical
quasi-geostrophic equation}

\author{Peter Constantin
\\Department of Mathematics \\University of Chicago
\\ 5734 S. University Avenue\\Chicago, IL 60637 \\ \\E-mail: const@cs.uchicago.edu
\\  \\ \\Jiahong Wu
\\Department of Mathematics\\ Oklahoma State University\\ Stillwater, OK 74078\\ \\
E-mail: jiahong@math.okstate.edu}
\date{}
\maketitle

\newtheorem{thm}{Theorem}[section]
\newtheorem{thmA}[thm]{Theorem A}
\newtheorem{cor}[thm]{Corollary}
\newtheorem{prop}[thm]{Proposition}
\newtheorem{define}[thm]{Definition}
\newtheorem{rem}[thm]{Remark}
\newtheorem{example}[thm]{Example}
\newtheorem{lemma}[thm]{Lemma}
\def\theequation{\thesection.\arabic{equation}}

\vspace{.25in} \noindent {\bf Abstract}. We present a
regularity result for weak solutions of the 2D quasi-geostrophic
equation with supercritical ($\alpha< 1/2$) dissipation $(-\Delta)^\alpha$ :
If a Leray-Hopf weak solution is H\"{o}lder
continuous $\theta\in C^\delta({\mathbb R}^2)$ with $\delta>1-2\alpha$
on the time interval $[t_0, t]$, then it is
actually a classical solution on $(t_0,t]$.

\vspace{.3in} \noindent AMS (MOS) Numbers: 76D03, 35Q35

\vspace{.2in} \noindent Keywords: the 2D quasi-geostrophic equation,
supercritical dissipation, regularity, weak solutions.

\newpage
\section{Introduction}
\setcounter{equation}{0} \label{sec:1}

We discuss the surface 2D quasi-geostrophic (QG) equation
\begin{equation}\label{qgse}
\partial_t  \theta +u\cdot\nabla\theta + \kappa
(-\Delta)^{\alpha}\theta=0,\quad x\in {\Bbb R}^2, \,\, t>0,
\end{equation}
where $\alpha> 0$ and $\kappa\ge 0$ are parameters, and the 2D
velocity field $u=(u_1, u_2)$ is determined from $\theta$ by the
stream function $\psi$ via the auxiliary relations
\begin{equation}\label{uth}
(u_1,u_2)=\left(-\partial_{x_2}\psi, \partial_{x_1}\psi\right),
\qquad {(-\Delta)^{\frac{1}{2}}}\psi=-\theta.
\end{equation}
Using the notation $\Lambda\equiv (-\Delta)^{\frac{1}{2}}$ and
$\nabla^\perp \equiv (\partial_{x_2}, -\partial_{x_1})$, the
relations in (\ref{uth}) can be combined into
\begin{equation}\label{uthr}
u = \nabla^\perp\,\Lambda^{-1}\, \theta =(-{\cal R}_2\theta, {\cal
R}_1\theta),
\end{equation}
where ${\cal R}_1$ and ${\cal R}_2$ are the usual Riesz transforms
in ${\Bbb R}^2$. The 2D QG equation with $\kappa>0$ and
$\alpha=\frac{1}{2}$ arises in geophysical studies of strongly
rotating fluids (see \cite{Con},\cite{Pe} and references therein)
while the inviscid QG equation ((\ref{qgse}) with $\kappa=0$) was
derived to model frontogenesis in meteorology, a formation of sharp
fronts between masses of hot and cold air (see
\cite{CMT},\cite{HPGS},\cite{Pe}).

\vspace{.07in} The problem at the center of the mathematical theory
concerning the 2-D QG equation is whether or not it has a global in
time smooth solution for any prescribed smooth initial data. In the
subcritical case $\alpha>\frac{1}{2}$, the dissipative QG equation
has been shown to possess a unique global smooth solution for every
sufficiently smooth initial data (see \cite{CW5},\cite{Res}). In
contrast, when $\alpha\le \frac12$, the issue of global existence
and uniqueness is more difficult and has still unanswered
aspects. Recently this problem
has attracted a significant amount of research
(\cite{CV},\cite{Cha},\cite{ChL},\cite{CMZ},\cite{Con},\cite{CCW},\cite{CC},
\cite{Ju},\cite{Ju2},\cite{KNV},\cite{MarLr},\cite{Sch},\cite{Sch2},\cite{Wu02},\cite{Wu3},
\cite{Wu4},\cite{Wu41},\cite{Wu77}). In Constantin, C\'{o}rdoba and
Wu \cite{CCW}, we proved in the critical case ($\alpha=\frac12$) the
global existence and uniqueness of classical solutions corresponding
to any initial data with $L^\infty$-norm comparable to or less than
the diffusion coefficient $\kappa$. In a recently posted preprint in
arXiv \cite{KNV}, Kiselev, Nazarov and Volberg proved that
smooth global solutions exist for any $C^\infty$ periodic initial
data, by removing the $L^\infty$-smallness assumption on the initial
data of \cite{CCW}. Caffarelli and Vasseur (arXiv reference
\cite{CV}) establish the global regularity of the Leray-Hopf type
weak solutions (in $L^\infty((0,\infty);L^2)\cap
L^2((0,\infty);\mathring{H}^{1/2})$) of the critical 2D QG equation with
$\alpha=\frac{1}{2}$ in general ${\Bbb R}^n$.

\vspace{.07in} In this paper we present a regularity result of weak
solutions of the dissipative QG equation with $\alpha<\frac12$ (the
supercritical case). The result asserts that if a Leray-Hopf weak solution
$\theta$  of (\ref{qgse}) is in the H\"{o}lder class $C^\delta$ with
$\delta>1-2\alpha$ on the time interval $[t_0,t]$, then it is
actually a classical solution on $(t_0,t]$. The proof involves
representing the functions in H\"{o}lder space in terms of the
Littlewood-Paley decomposition and using Besov space techniques. When
$\theta$ is in $C^\delta$, it also belongs to the Besov space
$\mathring{B}^{\delta(1-2/p)}_{p,\infty}$ for any $p\ge 2$. By
taking $p$ sufficiently large, we have $\theta\in C^{\delta_1}\cap
\mathring{B}^{\delta_1}_{p,\infty}$ for $\delta_1>1-2\alpha$. The
idea is to show that $\theta\in C^{\delta_2} \cap
\mathring{B}^{\delta_2}_{p,\infty}$ with $\delta_2>\delta_1$.
Through iteration, we establish that $\theta\in C^\gamma$ with
$\gamma>1$. Then $\theta$ becomes a classical solution.

\vspace{.1in} The results of this paper can be easily extended to a
more general form of the quasi-geostrophic equation in which $x\in
{\Bbb R}^n$ and $u$ is a divergence-free vector field determined by
$\theta$ through a singular integral operator.

\vspace{.1in} The rest of this paper is divided into two sections.
Section \ref{sec:2} provides the definition of Besov spaces and
necessary tools. Section \ref{sec:3} states and proves the main
result.

\vspace{.3in}
\section{Besov spaces and related tools}
\setcounter{equation}{0} \label{sec:2}

This section provides the definition of Besov spaces and several
related tools. We start with a some notation. Denote by ${\cal
S}({\Bbb R}^n)$ the usual Schwarz class and ${\cal S}'({\Bbb R}^n)$
the space of tempered distributions. $\widehat{f}$ denotes the
Fourier transform of $f$, namely
$$
\widehat{f}(\xi) = \int_{{\Bbb R}^n} e^{-ix\cdot\xi} \, f(x)\, dx.
$$
The fractional Laplacian $(-\Delta)^\alpha$ can be defined through
the Fourier transform
$$
\widehat{(-\Delta)^\alpha f} = |\xi|^{2\alpha} \, \widehat{f}(\xi).
$$
Let
$$
{\cal S}_0 = \left\{\phi\in {\cal S}, \int_{{\Bbb R}^n} \phi(x)
x^\gamma dx =0, \, |\gamma| = 0,1,2,\cdots \right\}.
$$
Its dual ${\cal S}_0'$  is given by
$$
{\cal S}_0' = {\cal S}'/{\cal S}_0^\perp = {\cal S}'/{\cal P},
$$
where ${\cal P}$ is the space of polynomials. In other words, two
distributions in ${\cal S}'$ are identified as the same in ${\cal
S}_0'$ if their difference is a polynomial.

\vspace{.1in} It is a classical result that there exists a dyadic
decomposition of ${\Bbb R}^n$, namely a sequence $\{\Phi_j\} \in
{\cal S}({\Bbb R}^n)$ such that
$$
\mbox{supp}\,\widehat{\Phi}_j \subset A_j,\quad \widehat{\Phi}_j
(\xi) = \widehat{\Phi}_0(2^{-j}\xi)\quad \mbox{or}\quad \Phi_j(x)
=2^{jn} \Phi_0(2^j x)
$$
and
$$
\sum_{k=-\infty}^\infty \widehat{\Phi}_k (\xi) =\left\{
\begin{array}{ll} 1 &\mbox{if $\xi \in {\Bbb R}^n\setminus\{0\}$},\\
0 &\mbox{if $\xi=0$},
\end{array}
\right.
$$
where
$$
A_j=\{\xi \in {\Bbb R}^n: \, 2^{j-1} < |\xi| < 2^{j+1}\}.
$$
As a consequence, for any $f\in {\cal S}_0'$,
\begin{equation}\label{id1}
\sum_{k=-\infty}^\infty \Phi_k \ast f = f.
\end{equation}
For notational convenience, set
\begin{equation}\label{d1}
\Delta_j f = \Phi_j \ast f, \quad j=0, \pm1, \pm2, \cdots.
\end{equation}
\begin{define}
For $s \in {\Bbb R}$ and $1\le p,q \le \infty$, the homogeneous
Besov space $\mathring{B}_{p,q}^s$ is defined by
$$
\mathring{B}_{p,q}^s = \left\{f\in {\cal
S}_0':\,\,\|f\|_{\mathring{B}_{p,q}^s}<\infty \right\},
$$
where
$$
\|f\|_{\mathring{B}_{p,q}^s} =\left\{
\begin{array}{ll}
\displaystyle \Big(\sum_{j}\, \Big(2^{js} \|\Delta_j
f\|_{L^p}\Big)^q \Big)^{1/q}  \quad  &\mbox{for $q<\infty$}, \\
\displaystyle \sup_j 2^{js} \, \|\Delta_j f\|_{L^p} \quad &\mbox{for
$q=\infty$}.
\end{array}
\right.
$$
\end{define}
For $\Delta_j$ defined in (\ref{d1}) and $S_j\equiv \sum_{k<j}
\Delta_k$,
$$
\Delta_j \Delta_k = 0\quad \mbox{if $|j-k|\ge 2$} \quad
\mbox{and}\quad \Delta_j(S_{k-1}f\,\Delta_k f) =0\quad\mbox{if
$|j-k|\ge 3$}.
$$
The following proposition lists a few simple facts that we will use
in the subsequent section.
\begin{prop}\label{embed}
Assume that $s\in {\Bbb R}$ and $p,q\in [1,\infty]$.
\begin{enumerate}
\item[1)] If $1\le q_1 \le q_2 \le \infty$, then $\mathring{B}^{s}_{p,q_1} \subset
\mathring{B}^s_{p,q_2}$.
\item[2)] ({\bf Besov embedding}) If $1\le p_1\le p_2\le \infty$ and
$s_1= s_2 + n(\frac{1}{p_1}-\frac{1}{p_2})$, then
$\mathring{B}^{s_1}_{p_1,q}({\Bbb R}^n) \subset
\mathring{B}^{s_2}_{p_2,q}({\Bbb R}^n)$.
\item[3)] If $1<p<\infty$, then
$$
\mathring{B}^s_{p,\min(p,2)} \subset \mathring{W}^{s,p} \subset
\mathring{B}^s_{p,\max(p,2)},
$$
where $\mathring{W}^{s,p}$ denotes a standard homogeneous Sobolev
space.
\end{enumerate}
\end{prop}

We will need a Bernstein type inequality for fractional derivatives.
\begin{prop}\label{bern}
Let $\alpha\ge0$. Let $1\le p\le q\le \infty$.
\begin{enumerate}
\item[1)] If $f$ satisfies
$$
\mbox{supp}\, \widehat{f} \subset \{\xi\in {\Bbb R}^n: \,\, |\xi|
\le K 2^j \},
$$
for some integer $j$ and a constant $K>0$, then
$$
\|(-\Delta)^\alpha f\|_{L^q({\Bbb R}^n)} \le C_1\, 2^{2\alpha j +
jn(\frac{1}{p}-\frac{1}{q})} \|f\|_{L^p({\Bbb R}^n)}.
$$
\item[2)] If $f$ satisfies
\begin{equation}\label{spp}
\mbox{supp}\, \widehat{f} \subset \{\xi\in {\Bbb R}^n: \,\, K_12^j
\le |\xi| \le K_2 2^j \}
\end{equation}
for some integer $j$ and constants $0<K_1\le K_2$, then
$$
C_1\, 2^{2\alpha j} \|f\|_{L^q({\Bbb R}^n)} \le \|(-\Delta)^\alpha
f\|_{L^q({\Bbb R}^n)} \le C_2\, 2^{2\alpha j +
jn(\frac{1}{p}-\frac{1}{q})} \|f\|_{L^p({\Bbb R}^n)},
$$
where $C_1$ and $C_2$ are constants depending on $\alpha,p$ and $q$
only.
\end{enumerate}
\end{prop}

The following proposition provides a lower bound for an integral
that originates from the dissipative term in the process of $L^p$
estimates (see \cite{Wu3},\cite{CMZ}).
\begin{prop}\label{lowerb}
Assume either $\alpha\ge 0$ and $p=2$ or $0\le \alpha\le 1$ and
$2<p<\infty$. Let $j$ be an integer and $f\in {\cal S}'$. Then
$$
\int_{{\Bbb R}^n} |\Delta_j f|^{p-2}\Delta_j f  \,\Lambda^{2\alpha}
\Delta_j f \,dx\ge C\, 2^{2\alpha j} \|\Delta_j f\|_{L^p}^p
$$
for some constant $C$ depending on $n$, $\alpha$ and $p$.
\end{prop}

\vspace{.3in}
\section{The main theorem and its proof}
\setcounter{equation}{0} \label{sec:3}

\begin{thm}
Let $\theta$ be a Leray-Hopf weak solution of (\ref{qgse}), namely
\begin{equation}\label{LH}
\theta\in L^\infty([0,\infty); L^2({\Bbb R}^2))\cap L^2([0,\infty);
\mathring{H}^\alpha({\Bbb R}^2)).
\end{equation}
Let $\delta>1-2\alpha$ and let $0<t_0<t<\infty$. If
\begin{equation}\label{holder}
\theta\in L^\infty([t_0,t]; C^\delta({\Bbb R}^2)),
\end{equation}
then
$$
\theta\in C^\infty((t_0,t]\times {\Bbb R}^2).
$$
\end{thm}
{\it Proof}.\quad First, we notice that (\ref{LH}) and
(\ref{holder}) imply that
$$
\theta\in L^\infty([t_0,t]; \mathring{B}^{\delta_1}_{p,\infty}({\Bbb
R}^2)),
$$
for any $p\ge 2$ and $\delta_1=\delta(1-\frac{2}{p})$. In fact, for
any $\tau\in [t_0,t]$,
\begin{eqnarray}
\|\theta(\cdot,\tau)\|_{\mathring{B}^{\delta_1}_{p,\infty}} &=&
\sup_{j}2^{\delta_1 j}\|\Delta_j \theta\|_{L^p} \nonumber \\
&\le& \sup_{j}2^{\delta_1 j}
\|\Delta_j\theta\|_{L^\infty}^{1-\frac{2}{p}} \,
\|\Delta_j\theta\|_{L^2}^{\frac2{p}}\nonumber \\
&\le& \|\theta(\cdot,\tau)\|_{C^\delta}^{1-\frac{2}{p}}\,
\|\theta(\cdot,\tau)\|_{L^2}^{\frac2{p}}.\nonumber
\end{eqnarray}
Since $\delta>1-2\alpha$, we have $\delta_1>1-2\alpha$ when
$$
p> p_0 \equiv \frac{2\delta}{\delta-(1-2\alpha)}.
$$

Next, we show that
$$
\theta\in L^\infty([t_0,t]; \mathring{B}^{\delta_1}_{p,\infty}\cap
C^{\delta_1})
$$
implies
$$
\theta(\cdot,t) \in \mathring{B}^{\delta_2}_{p,\infty}\cap
C^{\delta_2}
$$
for some $\delta_2>\delta_1$ to be specified. Let $j$ be an integer.
Applying $\Delta_j$ to (\ref{qgse}), we get
\begin{equation}\label{mot}
\partial_t \Delta_j \theta
+ \kappa \Lambda^{2\alpha} \Delta_j \theta =-\Delta_j (u\cdot \nabla
\theta).
\end{equation}
By Bony's notion of paraproduct,
\begin{eqnarray}
\Delta_j (u\cdot\nabla \theta) &=& \sum_{|j-k|\le 2} \Delta_j
(S_{k-1}u \cdot \nabla \Delta_k \theta) + \sum_{|j-k|\le 2} \Delta_j
(\Delta_k u \cdot \nabla S_{k-1} \theta) \nonumber \\ && +
\sum_{k\ge j-1}\sum_{|k-l|\le 1} \Delta_j (\Delta_k u \cdot\nabla
\Delta_l \theta). \label{para}
\end{eqnarray}
Multiplying  (\ref{mot}) by $p|\Delta_j \theta|^{p-2} \Delta_j
\theta$, integrating with respect to $x$, and applying the lower
bound
$$
\int_{{\Bbb R}^d} |\Delta_j f|^{p-2}\Delta_j f  \,\Lambda^{2\alpha}
\Delta_j f \,dx\ge C\, 2^{2\alpha j} \|\Delta_j f\|_{L^p}^p
$$
of Proposition \ref{lowerb}, we obtain
\begin{equation}\label{pri}
\frac{d}{dt} \|\Delta_j \theta\|^p_{L^p} + C\kappa 2^{2\alpha j}
\|\Delta_j \theta\|^p_{L^p}\le I_1 +I_2+I_3,
\end{equation}
where $I_1$, $I_2$ and $I_3$ are given by
\begin{eqnarray}
I_{1} &=& -p  \sum_{|j-k|\le 2} \int |\Delta_j \theta|^{p-2}
\Delta_j \theta\cdot \Delta_j (S_{k-1}u \cdot \nabla \Delta_k
\theta)\,dx, \nonumber\\
I_{2} &=& -p \sum_{|j-k|\le 2}  \int |\Delta_j \theta|^{p-2}
\Delta_j \theta\cdot \Delta_j (\Delta_k u \cdot \nabla S_{k-1}
\theta)\,dx, \nonumber\\
I_{3} &=& -p  \sum_{k\ge j-1}  \int |\Delta_j \theta|^{p-2} \Delta_j
\theta\cdot\sum_{|k-l|\le 1} \Delta_j (\Delta_k u\cdot\nabla
\Delta_l \theta)\,dx.\nonumber
\end{eqnarray}
We first bound $I_2$. By H\"{o}lder's inequality
$$
I_{2} \le C  \|\Delta_j \theta\|_{L^p}^{p-1} \sum_{|j-k|\le 2}
\|\Delta_k u\|_{L^p} \|\nabla S_{k-1}\theta\|_{L^\infty}.
$$
Applying Bernstein's inequality, we obtain
\begin{eqnarray}
I_{2} &\le& C \|\Delta_j \theta\|_{L^p}^{p-1} \sum_{|j-k|\le 2}
\|\Delta_k u\|_{L^p}
\sum_{m \le k-1} 2^m \|\Delta_m \theta\|_{L^\infty} \nonumber\\
&\le& C \|\Delta_j \theta\|_{L^p}^{p-1} \sum_{|j-k|\le 2} \|\Delta_k
u\|_{L^p} 2^{(1-\delta_1)k} \sum_{m\le k-1} 2^{(m-k)(1-\delta_1)}
2^{m\delta_1} \|\Delta_m \theta\|_{L^\infty}.\nonumber
\end{eqnarray}
Thus, for $1-\delta_1>0$, we have
$$
I_2 \le C \|\Delta_j \theta\|_{L^p}^{p-1}
\|\theta\|_{C^{\delta_1}}\sum_{|j-k|\le 2}\|\Delta_k u\|_{L^p}
2^{(1-\delta_1)k} .
$$
We now estimate $I_1$. The standard idea is to decompose it into
three terms: one with commutator, one that becomes zero due to the
divergence-free condition and the rest. That is, we rewrite $I_{1}$
as
\begin{eqnarray}
I_{1} &=& -p  \sum_{|j-k|\le 2}  \int |\Delta_j \theta|^{p-2}
\Delta_j \theta\cdot [\Delta_j, S_{k-1}u\cdot\nabla]\Delta_k
\theta\,dx \nonumber \\
&& -p  \int |\Delta_j \theta|^{p-2}
\Delta_j \theta\cdot (S_ju \cdot \nabla\Delta_j \theta)\,dx \nonumber \\
&& -p  \sum_{|j-k|\le 2}  \int |\Delta_j \theta|^{p-2} \Delta_j
\theta\cdot (S_{k-1}u -S_j u) \cdot\nabla \Delta_j\Delta_k
\theta\,dx
\nonumber \\
&=& I_{11} + I_{12} + I_{13},\nonumber
\end{eqnarray}
where we have used the simple fact that $\sum_{|k-j|\le 2}
\Delta_k\Delta_j\theta =\Delta_j\theta$, and the brackets $[\,]$
represent the commutator, namely
$$
[\Delta_j, S_{k-1}u\cdot\nabla]\Delta_k \theta
=\Delta_j(S_{k-1}u\cdot\nabla \Delta_k \theta)-
S_{k-1}u\cdot\nabla\Delta_j\Delta_k \theta.
$$
Since $u$ is divergence free, $I_{12}$ becomes zero. $I_{12}$ can
also be handled without resort to the divergence-free condition. In
fact, integrating by parts in $I_{12}$ yields
$$
I_{12} =\int |\Delta_j \theta|^p \,\nabla\cdot S_ju\, dx \le
\|\Delta_j \theta\|^p_{L^p} \|\nabla\cdot S_ju\|_{L^\infty}.
$$
By Bernstein's inequality,
\begin{eqnarray}
|I_{12}| &\le& \|\Delta_j \theta\|^p_{L^p} \sum_{m\le j-1} 2^m
\|\Delta_m u\|_{L^\infty} \nonumber \\
&=&\|\Delta_j \theta\|^p_{L^p} 2^{(1-\delta_1)j}\sum_{m\le j-1}
2^{(1-\delta_1)(m-j)}\,2^{m\delta_1}\|\Delta_m u\|_{L^\infty}.
\nonumber
\end{eqnarray}
For $1-\delta_1>0$,
$$
|I_{12}| \le C\, \|\Delta_j \theta\|^p_{L^p} 2^{(1-\delta_1)j}
\|u\|_{C^{\delta_1}} \le C\,\|\Delta_j
\theta\|^{p-1}_{L^p}\,2^{(1-2\delta_1)j}
\,\|\theta\|_{\mathring{B}^{\delta_1}_{p,\infty}}\|u\|_{C^{\delta_1}}.
$$
We now bound $I_{11}$ and $I_{13}$. By H\"{o}lder's inequality,
$$
|I_{11}|  \le p\|\Delta_j\theta\|_{L^p}^{p-1} \sum_{|j-k|\le
2}\|[\Delta_j, S_{k-1}u\cdot\nabla]\Delta_k \theta\|_{L^p}.
$$
To bound the the commutator, we have by the definition of $\Delta_j$
$$
[\Delta_j, S_{k-1}u\cdot\nabla]\Delta_k \theta =\int
\Phi_j(x-y)\left(S_{k-1}(u)(x)-S_{k-1}(u)(y)\right)\cdot\nabla
\Delta_k \theta(y)\,dy.
$$
Using the fact that $\theta\in C^{\delta_1}$ and thus
$$
\|S_{k-1}(u)(x)-S_{k-1}(u)(y)\|_{L^\infty} \le
\|u\|_{C^{\delta_1}}\,|x-y|^{\delta_1},
$$
we obtain
$$
\|[\Delta_j, S_{k-1}u\cdot\nabla]\Delta_k \theta\|_{L^p} \le
2^{-\delta_1 j}\,\|u\|_{C^{\delta_1}} 2^k\|\Delta_k \theta\|_{L^p}.
$$
Therefore,
$$
|I_{11}|  \le C p\,\|\Delta_j\theta\|_{L^p}^{p-1}\,2^{-\delta_1
j}\,\|u\|_{C^{\delta_1}}\sum_{|j-k|\le
2}2^k\|\Delta_k\theta\|_{L^p}.
$$
The estimate for $I_{13}$ is straightforward. By H\"{o}lder's
inequality,
\begin{eqnarray}
|I_{13}| &\le& p\|\Delta_j\theta\|_{L^p}^{p-1} \sum_{|j-k|\le 2}
\|S_{k-1}u -S_j u\|_{L^p} \|\nabla \Delta_j \theta\|_{L^\infty}
\nonumber \\ &\le& C
p\,\|\Delta_j\theta\|_{L^p}^{p-1}\,2^{(1-\delta_1)j}\,\|\theta\|_{C^{\delta_1}}\,
\sum_{|j-k|\le 2}\|\Delta_k u\|_{L^p}. \nonumber
\end{eqnarray}
We now bound $I_3$. By
H\"{o}lder's inequality and Bernstein's inequality,
\begin{eqnarray}
|I_3| &\le&  p \|\Delta_j\theta\|_{L^p}^{p-1} \,\|\Delta_j
\nabla\cdot\Big(\sum_{k\ge j-1} \sum_{|l-k|\le 1}\Delta_l u
\,\Delta_k
\theta\Big)\|_{L^p} \nonumber \\
&\le&  p \|\Delta_j\theta\|_{L^p}^{p-1} \,2^j \|u\|_{C^{\delta_1}}
\sum_{k\ge j-1} 2^{-\delta_1 k} \|\Delta_k\theta\|_{L^p}.
\end{eqnarray}
Inserting the estimates for $I_1$, $I_2$ and $I_3$ in (\ref{pri})
and eliminating $p\|\Delta_j\theta\|_{L^p}^{p-1}$ from both sides,
we get
\begin{eqnarray}
\frac{d}{dt} \|\Delta_j \theta\|_{L^p} + C\kappa 2^{2\alpha j}
\|\Delta_j \theta\|_{L^p} &\le&  C\,2^{(1-2\delta_1)j}
\,\|\theta\|_{\mathring{B}^{\delta_1}_{p,\infty}}\|u\|_{C^{\delta_1}}
\nonumber \\
&&+ C 2^{-\delta_1 j}\,\|u\|_{C^{\delta_1}}\sum_{|j-k|\le
2}2^k\|\Delta_k\theta\|_{L^p}\nonumber \\
&& + C\,\|\theta\|_{C^{\delta_1}}\sum_{|j-k|\le 2}\|\Delta_k
u\|_{L^p}
2^{(1-\delta_1)k} \nonumber \\
&&+C\,2^{(1-\delta_1)j}\,\|\theta\|_{C^{\delta_1}}\, \sum_{|j-k|\le
2}\|\Delta_k u\|_{L^p}\nonumber \\
&& + C\, 2^j
\|u\|_{C^{\delta_1}}\sum_{k\ge j-1} 2^{-\delta_1 k} \|\Delta_k
\theta\|_{L^p}. \label{good}
\end{eqnarray}
The terms on the right can be further bounded as follows.
\begin{eqnarray}
C 2^{-\delta_1 j}\,\|u\|_{C^{\delta_1}}\sum_{|j-k|\le
2}2^k\|\Delta_k\theta\|_{L^p} &=& C 2^{(1-2\delta_1)
j}\,\,\|u\|_{C^{\delta_1}}\sum_{|j-k|\le 2} 2^{\delta_1
k}\|\Delta_k\theta\|_{L^p}\,2^{(k-j)(1-\delta_1)}\nonumber \\
&\le& C\, 2^{(1-2\delta_1)j} \|u\|_{C^{\delta_1}}\,
\|\theta\|_{\mathring{B}^{\delta_1}_{p,\infty}}, \nonumber
\end{eqnarray}
\begin{eqnarray}
C\,\|\theta\|_{C^{\delta_1}}\sum_{|j-k|\le 2}\|\Delta_k u\|_{L^p}
2^{(1-\delta_1)k} &=& C 2^{(1-2\delta_1)
j}\,\,\|\theta\|_{C^{\delta_1}}\sum_{|j-k|\le 2} 2^{\delta_1
k}\|\Delta_k u\|_{L^p}\,2^{(k-j)(1-2\delta_1)}\nonumber \\
&\le& C\, 2^{(1-2\delta_1)j} \|\theta\|_{C^{\delta_1}}\,
\|u\|_{\mathring{B}^{\delta_1}_{p,\infty}}, \nonumber
\end{eqnarray}
\begin{eqnarray}
C\,2^{(1-\delta_1)j}\,\|\theta\|_{C^{\delta_1}}\, \sum_{|j-k|\le
2}\|\Delta_k u\|_{L^p} &=& C\,2^{(1-2\delta_1)
j}\,\|\theta\|_{C^{\delta_1}}\sum_{|j-k|\le 2}2^{\delta_1
k}\|\Delta_k u\|_{L^p}\,\,2^{(j-k)\delta_1}\nonumber \\
&\le& C\, 2^{(1-2\delta_1)j} \|\theta\|_{C^{\delta_1}}\,
\|u\|_{\mathring{B}^{\delta_1}_{p,\infty}}\nonumber
\end{eqnarray}
and
\begin{eqnarray}
C\, 2^j \|u\|_{C^{\delta_1}}\sum_{k\ge j-1} 2^{-\delta_1 k}
\|\Delta_k \theta\|_{L^p} &=&  C\, 2^{(1-2\delta_1)j}
\|u\|_{C^{\delta_1}}\sum_{k\ge j-1} 2^{-2\delta_1 (k-j)} 2^{\delta_1
k} \|\Delta_k \theta\|_{L^p}
\nonumber \\
&\le & C\, 2^{(1-2\delta_1)j} \|u\|_{C^{\delta_1}}\,
\|\theta\|_{\mathring{B}^{\delta_1}_{p,\infty}}. \nonumber
\end{eqnarray}
We can write (\ref{good}) in the following integral form
\begin{eqnarray}
\|\Delta_j \theta(t)\|_{L^p} &\le&  e^{-C \kappa\,2^{2\alpha j}
(t-t_0)} \|\Delta_j \theta(t_0)\|_{L^p} \nonumber \\
&& +\,C\,\int_{t_0}^t e^{-C\kappa 2^{2\alpha
j}\,(t-s)}2^{(1-2\delta_1)j} (\|\theta\|_{C^{\delta_1}}\,
\|u\|_{\mathring{B}^{\delta_1}_{p,\infty}}+ \|u\|_{C^{\delta_1}}\,
\|\theta\|_{\mathring{B}^{\delta_1}_{p,\infty}})\,ds . \nonumber
\end{eqnarray}
Multiplying both sides by $2^{(2\alpha+2\delta_1-1)j}$ and taking
the supremum with respect to $j$, we get
\begin{eqnarray}
\|\theta(t)\|_{\mathring{B}_{p,\infty}^{2\delta_1+2\alpha-1}} &\le&
\sup_{j} \{e^{-C \kappa\,2^{2\alpha j} (t-t_0)}
2^{(\delta_1+2\alpha-1) j}\}\,
\|\theta(t_0)\|_{\mathring{B}_{p,\infty}^{\delta_1}} \nonumber
\\ && + C\kappa^{-1} \sup_{j}\{(1-e^{-C \kappa\,2^{2\alpha j}
(t-t_0)})\}
\max_{s\in[t_0,t]}\|\theta(s)\|_{\mathring{B}_{p,\infty}^{\delta_1}}
\|\theta(s)\|_{C^{\delta_1}}\nonumber
\end{eqnarray}
Here we have used the fact that
$$
\|u\|_{C^{\delta_1}} \le \|\theta\|_{C^{\delta_1}} \quad\mbox{and}
\quad \|u\|_{\mathring{B}_{p,\infty}^{\delta_1}} \le
\|\theta\|_{\mathring{B}_{p,\infty}^{\delta_1}}
$$
Therefore, we conclude that if
$$
\theta \in L^\infty([t_0,t]; \mathring{B}^{\delta_1}_{p,\infty}\cap
C^{\delta_1}),
$$
then
\begin{equation}\label{bp2}
\theta(\cdot, t) \in \mathring{B}^{2\delta_1+2\alpha-1}_{p,\infty}.
\end{equation}
Since $\delta_1>1-2\alpha$, we have $2\delta_1+2\alpha-1>\delta_1$
and thus gain regularity. In addition, according to the Besov
embedding of Proposition \ref{embed},
$$
\mathring{B}^{2\delta_1+2\alpha-1}_{p,\infty} \subset
\mathring{B}^{\delta_2}_{\infty,\infty},
$$
where
$$
\delta_2 = 2\delta_1+2\alpha-1 -\frac{2}{p} = \delta_1 +
\left(\delta_1-\left(1-2\alpha+\frac2{p}\right)\right).
$$
We have $\delta_2>\delta_1$ when
$$
p>p_1\equiv \frac2{\delta_1-(1-2\alpha)}.
$$
Noting that
$$
\mathring{B}^{\delta_2}_{\infty,\infty}\cap L^\infty =C^{\delta_2},
$$
we conclude that, for $p>\max\{p_0,p_1\}$,
$$
\theta(\cdot,t) \in \mathring{B}^{\delta_2}_{p,\infty} \cap
C^{\delta_2}
$$
for some $\delta_2>\delta_1$. The above process can then be iterated
with $\delta_1$ replaced by $\delta_2$. A finite number of
iterations allow us to obtain that
$$
\theta(\cdot,t) \in C^\gamma
$$
for some $\gamma>1$. The regularity in the spatial variable can then
be converted into regularity in time. We have thus established that
$\theta$ is a classical solution to the supercritical QG equation.
Higher regularity can be proved by well-known methods.

\vspace{.4in}\noindent {\bf Acknowledgment}: PC was partially
supported by NSF-DMS 0504213. JW thanks the Department of
Mathematics at the University of Chicago for its support
and hospitality.


\vspace{.3in}
\begin{thebibliography}{99}
\bibitem{CV} L. Caffarelli and A. Vasseur, Drift diffusion
equations with fractional diffusion and the quasi-geostrophic
equation, ArXiv: Math.AP/0608447 (2006).

\bibitem{Cha}D. Chae, On the regularity conditions for the
dissipative quasi-geostrophic equations, {\it SIAM J. Math. Anal.}
{\bf 37} (2006), 1649-1656.

\bibitem{ChL} D. Chae and J. Lee, Global well-posedness in the super-critical
dissipative quasi-geostrophic equations, {\it Commun. Math. Phys.}
{\bf 233} (2003), 297-311.

\bibitem{CMZ} Q. Chen, C. Miao and Z. Zhang, A new Bernstein
inequality and the 2D dissipative quasi-geostrophic equation, to
appear in {\it Commun. Math. Phys.}.


\bibitem{Con} P. Constantin, Euler equations, Navier-Stokes equations and turbulence.
Mathematical foundation of turbulent viscous flows,  1--43, Lecture
Notes in Math., 1871, Springer, Berlin, 2006.

\bibitem{CCW} P. Constantin, D. Cordoba and J. Wu, On the critical dissipative
quasi-geostrophic equation, {\em Indiana Univ. Math. J.} {\bf 50}
(2001), 97-107.

\bibitem{CMT} P. Constantin, A. Majda, and  E. Tabak,
Formation of strong fronts in the 2-D quasi-geostrophic thermal
active scalar, {\it Nonlinearity} {\bf 7}(1994), 1495-1533.

\bibitem{CW5} P. Constantin and J. Wu, Behavior of solutions of 2D
quasi-geostrophic equations, {\it SIAM J. Math. Anal.} {\bf 30}
(1999), 937-948.


\bibitem{CC} A. C\'{o}rdoba and D. C\'{o}rdoba, A maximum principle
applied to quasi-geostrophic equations, {\it Commun. Math. Phys.}
{\bf 249} (2004), 511-528.

\bibitem{HPGS} I. Held, R. Pierrehumbert, S. Garner, and K. Swanson,
Surface quasi-geostrophic dynamics, {\it J. Fluid Mech.} {\bf 282}
(1995), 1-20.

\bibitem{Ju} N. Ju, The maximum principle and the global attractor for the dissipative
2D quasi-geostrophic equations,  {\it Commun. Math. Phys. \bf 255}
(2005), 161-181.

\bibitem{Ju2} N. Ju, Global solutions to the two dimensional quasi-geostrophic
equation with critical or super-critical dissipation, {\it Math.
Ann.} {\bf 334} (2006), 627--642.

\bibitem{KNV} A. Kiselev, F. Nazarov and A. Volberg, Global
well-posedness for the critical 2D dissipative quasi-geostrophic
equation, ArXiv: Math.AP/0604185 (2006).

\bibitem{MarLr} F. Marchand and P.G. Lemari\'{e}-Rieusset, Solutions
auto-similaires non radiales pour l'\'{e}quation
quasi-g\'{e}ostrophique dissipative critique, {\it  C. R. Math.
Acad. Sci. Paris} {\bf  341}  (2005),  535--538.

\bibitem{Pe}  J. Pedlosky, ``Geophysical fluid dynamics", Springer, New York,
1987.

\bibitem{Res} S. Resnick, Dynamical problems in nonlinear advective
partial differential equations, Ph.D. thesis, University of Chicago,
1995.

\bibitem{Sch} M. Schonbek and T. Schonbek, Asymptotic behavior to
dissipative quasi-geostrophic flows, {\it SIAM J. Math. Anal.} {\bf
35} (2003), 357-375.

\bibitem{Sch2} M. Schonbek and T. Schonbek, Moments and lower
bounds in the far-field of solutions to quasi-geostrophic flows,
{\it Discrete Contin. Dyn. Syst.} {\bf 13} (2005), 1277-1304.

\bibitem{Wu02}J. Wu, The quasi-geostrophic equation and its two
regularizations, {\it Commun. Partial Differential Equations}\, {\bf
27} (2002), 1161-1181.

\bibitem{Wu3} J. Wu, Global solutions of the 2D dissipative quasi-geostrophic equation in
Besov spaces, {\it SIAM J. Math. Anal.}\,\, {\bf 36} (2004/2005),
1014-1030.

\bibitem{Wu4} J. Wu, The quasi-geostrophic equation with critical or supercritical
dissipation, {\it Nonlinearity} \,\,{\bf 18} (2005), 139-154.

\bibitem{Wu41} J. Wu, Solutions of the 2-D quasi-geostrophic equation in H\"{o}lder
spaces, {\it Nonlinear Analysis}\,\, {\bf 62} (2005), 579-594.

\bibitem{Wu77} J. Wu, Existence and uniqueness results for the 2-D dissipative
quasi-geostrophic equation, {\it Nonlinear Analysis}, in press.


\end{thebibliography}
\end{document}